\documentclass[reqno,12pt,a4paper]{amsart}

\usepackage{mathrsfs}
\usepackage{amssymb}
\usepackage{amsfonts}
\usepackage{latexsym}
\usepackage{amsthm}

\usepackage{graphicx}
\def\lb{\label}

\newcommand{\er}[1]{\textrm{(\ref{#1})}}

\begin{document}


\renewcommand{\theequation}{\arabic{section}.\arabic{equation}}
\theoremstyle{plain}
\newtheorem{theorem}{\bf Theorem}[section]
\newtheorem{lemma}[theorem]{\bf Lemma}
\newtheorem{corollary}[theorem]{\bf Corollary}
\newtheorem{proposition}[theorem]{\bf Proposition}
\newtheorem{definition}[theorem]{\bf Definition}
\newtheorem{remark}[theorem]{\it Remark}

\def\a{\alpha}  \def\cA{{\mathcal A}}     \def\bA{{\bf A}}  \def\mA{{\mathscr A}}
\def\b{\beta}   \def\cB{{\mathcal B}}     \def\bB{{\bf B}}  \def\mB{{\mathscr B}}
\def\g{\gamma}  \def\cC{{\mathcal C}}     \def\bC{{\bf C}}  \def\mC{{\mathscr C}}
\def\G{\Gamma}  \def\cD{{\mathcal D}}     \def\bD{{\bf D}}  \def\mD{{\mathscr D}}
\def\d{\delta}  \def\cE{{\mathcal E}}     \def\bE{{\bf E}}  \def\mE{{\mathscr E}}
\def\D{\Delta}  \def\cF{{\mathcal F}}     \def\bF{{\bf F}}  \def\mF{{\mathscr F}}
\def\c{\chi}    \def\cG{{\mathcal G}}     \def\bG{{\bf G}}  \def\mG{{\mathscr G}}
\def\z{\zeta}   \def\cH{{\mathcal H}}     \def\bH{{\bf H}}  \def\mH{{\mathscr H}}
\def\e{\eta}    \def\cI{{\mathcal I}}     \def\bI{{\bf I}}  \def\mI{{\mathscr I}}
\def\p{\psi}    \def\cJ{{\mathcal J}}     \def\bJ{{\bf J}}  \def\mJ{{\mathscr J}}
\def\vT{\Theta} \def\cK{{\mathcal K}}     \def\bK{{\bf K}}  \def\mK{{\mathscr K}}
\def\k{\kappa}  \def\cL{{\mathcal L}}     \def\bL{{\bf L}}  \def\mL{{\mathscr L}}
\def\l{\lambda} \def\cM{{\mathcal M}}     \def\bM{{\bf M}}  \def\mM{{\mathscr M}}
\def\L{\Lambda} \def\cN{{\mathcal N}}     \def\bN{{\bf N}}  \def\mN{{\mathscr N}}
\def\m{\mu}     \def\cO{{\mathcal O}}     \def\bO{{\bf O}}  \def\mO{{\mathscr O}}
\def\n{\nu}     \def\cP{{\mathcal P}}     \def\bP{{\bf P}}  \def\mP{{\mathscr P}}
\def\r{\rho}    \def\cQ{{\mathcal Q}}     \def\bQ{{\bf Q}}  \def\mQ{{\mathscr Q}}
\def\s{\sigma}  \def\cR{{\mathcal R}}     \def\bR{{\bf R}}  \def\mR{{\mathscr R}}
\def\S{\Sigma}  \def\cS{{\mathcal S}}     \def\bS{{\bf S}}  \def\mS{{\mathscr S}}
\def\t{\tau}    \def\cT{{\mathcal T}}     \def\bT{{\bf T}}  \def\mT{{\mathscr T}}
\def\f{\phi}    \def\cU{{\mathcal U}}     \def\bU{{\bf U}}  \def\mU{{\mathscr U}}
\def\F{\Phi}    \def\cV{{\mathcal V}}     \def\bV{{\bf V}}  \def\mV{{\mathscr V}}
\def\P{\Psi}    \def\cW{{\mathcal W}}     \def\bW{{\bf W}}  \def\mW{{\mathscr W}}
\def\o{\omega}  \def\cX{{\mathcal X}}     \def\bX{{\bf X}}  \def\mX{{\mathscr X}}
\def\x{\xi}     \def\cY{{\mathcal Y}}     \def\bY{{\bf Y}}  \def\mY{{\mathscr Y}}
\def\X{\Xi}     \def\cZ{{\mathcal Z}}     \def\bZ{{\bf Z}}  \def\mZ{{\mathscr Z}}
\def\O{\Omega}

\newcommand{\gA}{\mathfrak{A}}          \newcommand{\ga}{\mathfrak{a}}
\newcommand{\gB}{\mathfrak{B}}          \newcommand{\gb}{\mathfrak{b}}
\newcommand{\gC}{\mathfrak{C}}          \newcommand{\gc}{\mathfrak{c}}
\newcommand{\gD}{\mathfrak{D}}          \newcommand{\gd}{\mathfrak{d}}
\newcommand{\gE}{\mathfrak{E}}
\newcommand{\gF}{\mathfrak{F}}           \newcommand{\gf}{\mathfrak{f}}
\newcommand{\gG}{\mathfrak{G}}           
\newcommand{\gH}{\mathfrak{H}}           \newcommand{\gh}{\mathfrak{h}}
\newcommand{\gI}{\mathfrak{I}}           \newcommand{\gi}{\mathfrak{i}}
\newcommand{\gJ}{\mathfrak{J}}           \newcommand{\gj}{\mathfrak{j}}
\newcommand{\gK}{\mathfrak{K}}            \newcommand{\gk}{\mathfrak{k}}
\newcommand{\gL}{\mathfrak{L}}            \newcommand{\gl}{\mathfrak{l}}
\newcommand{\gM}{\mathfrak{M}}            \newcommand{\gm}{\mathfrak{m}}
\newcommand{\gN}{\mathfrak{N}}            \newcommand{\gn}{\mathfrak{n}}
\newcommand{\gO}{\mathfrak{O}}
\newcommand{\gP}{\mathfrak{P}}             \newcommand{\gp}{\mathfrak{p}}
\newcommand{\gQ}{\mathfrak{Q}}             \newcommand{\gq}{\mathfrak{q}}
\newcommand{\gR}{\mathfrak{R}}             \newcommand{\gr}{\mathfrak{r}}
\newcommand{\gS}{\mathfrak{S}}              \newcommand{\gs}{\mathfrak{s}}
\newcommand{\gT}{\mathfrak{T}}             \newcommand{\gt}{\mathfrak{t}}
\newcommand{\gU}{\mathfrak{U}}             \newcommand{\gu}{\mathfrak{u}}
\newcommand{\gV}{\mathfrak{V}}             \newcommand{\gv}{\mathfrak{v}}
\newcommand{\gW}{\mathfrak{W}}             \newcommand{\gw}{\mathfrak{w}}
\newcommand{\gX}{\mathfrak{X}}               \newcommand{\gx}{\mathfrak{x}}
\newcommand{\gY}{\mathfrak{Y}}              \newcommand{\gy}{\mathfrak{y}}
\newcommand{\gZ}{\mathfrak{Z}}             \newcommand{\gz}{\mathfrak{z}}

\def\ve{\varepsilon}   \def\vt{\vartheta}    \def\vp{\varphi}    \def\vk{\varkappa}

\def\Z{{\mathbb Z}}    \def\R{{\mathbb R}}   \def\C{{\mathbb C}}\def\K{{\mathbb K}}
\def\T{{\mathbb T}}    \def\N{{\mathbb N}}   \def\dD{{\mathbb D}}
\def\B{{\mathbb B}}


\def\la{\leftarrow}              \def\ra{\rightarrow}      \def\Ra{\Rightarrow}
\def\ua{\uparrow}                \def\da{\downarrow}
\def\lra{\leftrightarrow}        \def\Lra{\Leftrightarrow}


\def\lt{\biggl}                  \def\rt{\biggr}
\def\ol{\overline}               \def\wt{\widetilde}
\def\no{\noindent}


\let\ge\geqslant                 \let\le\leqslant
\def\lan{\langle}                \def\ran{\rangle}
\def\/{\over}                    \def\iy{\infty}
\def\sm{\setminus}               \def\es{\emptyset}
\def\ss{\subset}                 \def\ts{\times}
\def\pa{\partial}                \def\os{\oplus}
\def\om{\ominus}                 \def\ev{\equiv}
\def\iint{\int\!\!\!\int}        \def\iintt{\mathop{\int\!\!\int\!\!\dots\!\!\int}\limits}
\def\el2{\ell^{\,2}}             \def\1{1\!\!1}
\def\sh{\sharp}
\def\wh{\widehat}

\def\where{\mathop{\mathrm{where}}\nolimits}
\def\as{\mathop{\mathrm{as}}\nolimits}
\def\Area{\mathop{\mathrm{Area}}\nolimits}
\def\arg{\mathop{\mathrm{arg}}\nolimits}
\def\const{\mathop{\mathrm{const}}\nolimits}
\def\det{\mathop{\mathrm{det}}\nolimits}
\def\diag{\mathop{\mathrm{diag}}\nolimits}
\def\diam{\mathop{\mathrm{diam}}\nolimits}
\def\dim{\mathop{\mathrm{dim}}\nolimits}
\def\dist{\mathop{\mathrm{dist}}\nolimits}
\def\Im{\mathop{\mathrm{Im}}\nolimits}
\def\Iso{\mathop{\mathrm{Iso}}\nolimits}
\def\Ker{\mathop{\mathrm{Ker}}\nolimits}
\def\Lip{\mathop{\mathrm{Lip}}\nolimits}
\def\rank{\mathop{\mathrm{rank}}\limits}
\def\Ran{\mathop{\mathrm{Ran}}\nolimits}
\def\Re{\mathop{\mathrm{Re}}\nolimits}
\def\Res{\mathop{\mathrm{Res}}\nolimits}
\def\res{\mathop{\mathrm{res}}\limits}
\def\sign{\mathop{\mathrm{sign}}\nolimits}
\def\span{\mathop{\mathrm{span}}\nolimits}
\def\supp{\mathop{\mathrm{supp}}\nolimits}
\def\Tr{\mathop{\mathrm{Tr}}\nolimits}
\def\BBox{\hspace{1mm}\vrule height6pt width5.5pt depth0pt \hspace{6pt}}


\newcommand\nh[2]{\widehat{#1}\vphantom{#1}^{(#2)}}
\def\dia{\diamond}

\def\Oplus{\bigoplus\nolimits}



\def\qqq{\qquad}
\def\qq{\quad}
\let\ge\geqslant
\let\le\leqslant
\let\geq\geqslant
\let\leq\leqslant
\newcommand{\ca}{\begin{cases}}
\newcommand{\ac}{\end{cases}}
\newcommand{\ma}{\begin{pmatrix}}
\newcommand{\am}{\end{pmatrix}}
\renewcommand{\[}{\begin{equation}}
\renewcommand{\]}{\end{equation}}
\def\bu{\bullet}

\title[{Trace formulae for Schr\"odinger operators on lattice}]
{Trace formulae for Schr\"odinger operators on lattice}

\date{\today}

\author[Evgeny Korotyaev]{Evgeny L. Korotyaev}
\address{Saint-Petersburg State University,   Universitetskaya
nab. 7/9, St. Petersburg, 199034, Russia, e.korotyaev@spbu.ru, \
korotyaev@gmail.com }


\begin{abstract}

We consider  Schr\"odinger operators with complex decaying
potentials (in general, not from trace class) on the lattice. We
determine trace formulae and  estimate of eigenvalues and singular
measure in terms of potentials. The proof is based on estimates of
free resolvent and analysis of functions  from Hardy space.
\end{abstract}

\maketitle

\section{Introduction}

We consider Schr\"odinger operators $H=\D+V$ on the lattice $\Z^{d},
d\ge 3$, where $\D$ is the discrete Laplacian on $\ell^{2}({\Z}^d)$
given by
\[
\big(\D f\big)(n)=\frac{1}{2}\sum_{|n-m|=1} f_m, \qqq
n=(n_j)_1^d\in\Z^d,
\]
where $f=(f_n)_{n\in{\Z}^d} \in \ell^{2}({\Z}^d)$. We assume that
the potential $V$  satisfies
\begin{equation}
\label{V} (V f)(n) =V_nf_n,\qqq
V\in \ell^p(\Z^{d}),\qqq\ca 1\le p<{6\/5} & if \ \ d=3\\
1\le p<{4\/3} & if  \  \ d\ge 4 \ac.
\end{equation}
Here  $\ell^p(\Z^{d}), p\ge 1$ is the  Banach space of sequences
 $f=(f_n)_{n\in \Z^d}$ equipped with the
norm
$$
\begin{aligned}
\|f\|_{p}=\|f\|_{\ell^p(\Z^{d})} =
\begin{cases}   \sup_{n\in \Z^d}|f_n|,\qq &
\ p=\iy, \\
\big(\sum_{n\in \Z^d}|f_n|^p\big)^{1\/p},\qq & \ p\in [1,\iy).
\end{cases}
\end{aligned}
$$
It is well-known that the spectrum of the Laplacian $\D$ is
absolutely continuous and equals
$$
\s(\D)=\s_{\textup{ac}}(\D)=[-d,d].
$$
Note that if $V$ satisfies \eqref{V}, then $V$ is the
Gilbert-Schmidt  operator and thus the essential spectrum of the
Schr\"odinger operator $H$ is given by $
\s_{\textup{ess}}(H)=[-d,d]. $ The operator $H$ has $N\le \iy$
eigenvalues $\{\l_n, n=1,....,N\}$ outside the interval $[-d,d]$.

We define the disc $\dD_r\subset\C$ with the radius $r>0$   by
$$
\dD_r=\{z\in \C:|z|<r\},
$$
and abbreviate  $\dD=\dD_1$. Define {\bf the new spectral variable}
$z\in \dD$ by
$$
\l=\l(z)={d\/2}\rt(z+{1\/z}\rt)\in \L=\C\sm [-d,d] ,\qqq z\in \dD.
$$
 The function $\l(z)$ has  the following properties (see more in Section 3):

{\it $\bu$  The function $\l(z)$ is a conformal mapping from $\dD$
onto the spectral domain $\L$.

$\bu$  The function $\l(z)$ maps the point $z=0$ to the point
$\l=\iy$.

$\bu$  The inverse mapping $z(\cdot ): \L\to \dD$ is given by
$z={1\/d}\rt(\l-\sqrt{\l^2-d^2}\rt),\ \ \l\in \L$ defined by
asymptotics $z={d\/2\l}+{O(1)\/\l^3}$ as $|\l|\to \iy$. }

Define the Hardy space $ \mH_p=\mH_p(\dD), 0<p\le \iy$. Let $F$ be
analytic in $\dD$. We say $F$ belongs the Hardy space $ \mH_p$ if
$F$ satisfies $\|F\|_{\mH_p}<\iy$, where $\|F\|_{\mH_p}$ is given by
$$
\|F\|_{\mH_p}=\ca
\sup_{0<r<1}\rt({1\/2\pi}\int_\T |F(re^{i\vt})|^pd\vt\rt)^{1\/p} &
if  \qqq 0< p<\iy\\
 \sup_{z\in \dD}|F(z)| & if \qqq p=\iy\ac ,
$$
where $\T=\R/(2\pi \Z)$. Note that the definition of the Hardy space
$\mH_p$ involves all $r\in (0,1)$.

Let  $\cB$ denote the class of bounded operators. Let $\cB_1$ and
$\cB_2$ be the trace and the Hilbert-Schmidt classes equipped with
the norm $\|\cdot \|_{\cB_1}$ and $ \|\cdot \|_{\cB_2}$
correspondingly.

Introduce the free resolvent $ R_0(k)=(\D-\l)^{-1}, \ \l\in \L$. For
$V\in \ell^2(\Z^d)$ we define the regularized determinant $\cD(\l)$
in the cut domain $\L$ and the modified determinant  $D$ in the disc
$\dD$ by
\[
\lb{De2}
\begin{aligned}
& \cD(\l)=\det \rt[(I+VR_0(\l))e^{-VR_0(\l)}\rt],\qqq \l\in\L,
\\
& D(z)=\cD(\l(z)),\qq z\in\dD.
\end{aligned}
\]
The function $\cD$ is a suitable regularization of the non-defined
determinant $\det (I+VR_0(\l))$, see \cite{GK69}. The function $\cD$
is analytic in $\L$ and the function $D$ is analytic in the disc
$\dD$. It has $N\le \iy$ zeros (counted with multiplicity) $z_1,
z_2,...$ in the disc $\dD$. Note that $\l_j=\l(z_j)$ is en
eigenvalue of $H$ (counted with multiplicity).

\subsection{Complex potentials}
In this paper we combine  classical results about Hardy spaces and
estimates of the free resolvent from \cite{KM17}, this gives us new
trace formulae for discrete Scr\"odinger operators $H=\D+V$ on the 
lattice $\Z^d$, where the potential $V$ is complex and satisfies the condition 
\er{V}. We improve results from \cite{KL16}, where potentials are considered under the weaker
condition $|V|^{2\/3}\in \ell^1(\Z^d)$.

 Introduce the additional conformal mapping $\vk:\L\to \K=\C\sm
[id,-id]$ by
\[
\lb{dvk}
\begin{aligned}
\vk=\sqrt{\l^2-d^2},\qqq \l\in \L,
\\
 \vk=\l-{d^2\/2\l}+{O(1)\/\l^3}\qq \as \
|\l|\to \iy.
\end{aligned}
\]
Note that if $\l\in\L$ and $\l\to \l_0\in [-d,d]$, then we obtain
$|\Im \l|+|\Re \vk(\l)|\to 0$.

\begin{theorem}
\lb{T1} Let a potential $V$ satisfy \er{V}. Then the modified
determinant $D$ is analytic in  the disc  $\Bbb D$ and is H\"older
up to the boundary and satisfies
\begin{equation}
\label{D1}
\begin{aligned}
\|D\|_{\mH_\iy(\dD)}\le e^{C_*^2\|V\|_{p}^2/2},
\end{aligned}
\end{equation}
where the constant $C_*$ is defined in \er{Cpd}. It has $N\le \iy$
zeros $\{z_j\}_{j=1}^N$ in the disc $\dD$, such that
\[
\lb{D2}
\begin{aligned}
& 0<r_0=\inf |z_j|=|z_1|\le |z_2|\le ...\le |z_j|\le |z_{j+1}\le
|z_{j+2}|\le .... ,\\
 & \sum _{j=1}^N \rt[(1-|z_j|)+|\Im \l_j|+|\Re \vk(\l_j)|  \rt]<\iy.
\end{aligned}
\]
 Moreover, the function $\p(z)= \log D(z)$
is analytic in $\dD_{r_0}$ and has the following Tailor series (here
$a={2\/d}$)
\[
\label{D3}
\begin{aligned}
 \p(z)=-\p_2z^2-\p_3z^3-\p_4z^4 +......, \qqq \as \qq
|z|<r_0,\\
\p_2= {a^2\/2}\Tr V^2, \qq \p_3={a^3\/3}\Tr\,V^3, \qq
\p_4={a^4\/4}{\rm Tr}\,(V^4+2VH_0VH_0+4V^2H_0^2)-\p_2,.....
\end{aligned}
\]

\end{theorem}

\medskip

For the function $D$ we define the Blaschke product $B(z), z\in \dD$
by: $B=1$ if $N=0$ and
\[
\lb{B2}
\begin{aligned}
 B(z)=\prod_{j=1}^N {|z_j|\/z_j}{(z_j-z)\/(1-\ol z_j z)},\qqq
  if \qqq N\ge 1.
\end{aligned}
\]
It is well known that the Blaschke product $B(z), z\in \dD$ given by
\er{B2} converges absolutely for $\{|z|<1\}$ and  satisfies $B\in
\mH_\iy$ with $\|B\|_{\mH_\iy}\le 1$, since $D\in \mH_\iy$
\cite{Koo98}. The Blaschke product $B$ has the Taylor series at
$z=0$:
\[
\begin{aligned}
\lb{B6}
& \log  B(z)=B_0-B_1z-B_2z^2-... \qqq as \qqq z\to 0,\\
& B_0=\log  B(0)<0,\qqq B_1=\sum_{j=1}^N\rt({1\/z_j}-\ol z_j
\rt),..., \qqq B_n={1\/n}\sum_{j=1}^N\rt({1\/z_j^n}-\ol z_j^n
\rt),....
\end{aligned}
\]
where each $B_n$ satisfy $|B_n|\le {2\/r_0^n}\sum _{j=1}^N
(1-|z_j|)$.

\bigskip
\noindent We describe the canonical representation of the
determinant $D(z)=\cD(\l(z)), z\in\dD$.

\begin{corollary}
\lb{T2} Let a potential $V$ satisfy \er{V}. Then there exists a
singular measure $\s\ge 0$ on $[-\pi,\pi]$, such that the
determinant $D$ has a canonical factorization for all $|z|<1$ given
by
\[
\lb{cfD}
\begin{aligned}
& D(z)=B(z)e^{-K_\s (z)}e^{K_D(z)},\\
& K_\s(z)={1\/2\pi}\int_{-\pi}^{\pi}{e^{it}+z\/e^{it}-z}d\s(t),\\
& K_D(z)= {1\/2\pi}\int_{-\pi}^{\pi}{e^{it}+z\/e^{it}-z}\log
|D(e^{it})|dt,
 \end{aligned}
\]
 where $\log |D(e^{it}) |\in L^1(-\pi,\pi)$ and
the measure $\s$ satisfies
\[
\supp \s\ss \{t\in [-\pi,\pi]: D(e^{it})=0\}.
\]
\end{corollary}

\noindent {\bf Remarks.}
1) For the canonical factorisation of analytic functions see, for
example, \cite{Koo98}.

\smallskip
2) Note that for the inner function $D_{in}(z)$ defined by
$D_{in}(z)= B(z) e^{-K_\sigma(z)}$, we have  $| D_{in}(z)|\le 1$,
since $d\s\ge 0$ and  $\Re {e^{it}+z\/e^{it}-z}\ge 0$ for all
$(t,z)\in \T\ts\dD$.

\smallskip
3) The function $D_B={D\/B}$  has no zeros in the disk $\dD$ and
satisfies
$$
\log D_B(z)={1\/2\pi}\int_{-\pi}^{\pi}{e^{it}+z\/e^{it}-z}d\m(t),
\qqq z\in \dD,
$$
where the measure $\m$ equals
$$
 d\m(t)=\log
|D(e^{it})|dt-d\s(t).
$$

\smallskip
\begin{theorem}
\lb{T3} {\bf (Trace formulae.)}  Let $V$ satisfy \er{V}. Then the
following identities hold:
\[
\lb{t1} {\s(\T)\/2\pi}-B_0={1\/2\pi}\int_{-\pi}^{\pi}\log
|D(e^{it})|dt\ge 0,
\]
\[
\lb{t2} B_1=\sum_{j=1}^N\rt({1\/z_j}-\ol z_j \rt)={1\/\pi}\int_\T
e^{-it}d\m(t),
\]
\[
\lb{t4}  \sum_{j=1}^N\rt({1\/z_j^2}-\ol z_j^2 \rt)={2\/d^2}\Tr \,
V^2+{1\/\pi}\int_\T e^{-i2t}d\m(t),
\]
\[
\lb{t3} B_n=\p_n+{1\/\pi}\int_\T e^{-int}d\m(t),\qqq n=2,3,....
\]
where $B_0=\log B(0)=\log \left(\prod_{j=1}^N |z_j|\right)<0$ and
$B_n$ are given by \er{B6}. In particular,
\[
\lb{zxj} \sum_{j=1}^N\rt(\Re \vk_j+i \Im \l_j\rt)={d\/2\pi}\int_\T
e^{-it}d\m(t),
\]
\[
\lb{zx}
\begin{aligned}
  \sum_{j=1}^N\rt[(\vk_{j1}\l_{j1}-\l_{j2}\vk_{j2})+i
  (\l_{j1}\l_{j2}+\vk_{j1}\vk_{j2})\rt]={1\/2}\Tr \,
V^2+{d^2\/4\pi}\int_\T e^{-2it}d\m(t),\\
\end{aligned}
\]
where $\vk_j=\sqrt{\l_j^2-d^2}=\vk_{j1}+i\vk_{j2}$ and
$\l_j=\l_{j1}+i\l_{j2}$.
\end{theorem}


We describe estimates of eigenvalues in terms of potentials.

\smallskip
\begin{theorem}
\lb{T4} Let  $V$ satisfy \er{V}. Then we have the following
estimates:
\[
\lb{eiV} \sum (1-|z_j|)\le -B_0\le
{C_*^2\/2}\|V\|_p^2-{\s(\T)\/2\pi},
\]
\[
\lb{eiV1}
\begin{aligned}
 |\sum_{j=1}^N\big(\Re \vk_j+i \Im \l_j\big)|\le d C_*^2\|V\|_p^2,
\end{aligned}
\]
\[
\lb{em}
\begin{aligned}
  {1\/2\pi}\big|\int_\T
e^{-int}d\m(t)\big|\le C_*^2\|V\|_p^2, \qqq \forall \ n\in \Z,
\end{aligned}
\]
and in particular,
\[
\lb{eiV2}
\begin{aligned}
 \sum_{j=1}^N\Im \l_j\le dC_*^2\|V\|_p^2,\qqq if \qqq \Im V\ge 0,
 \\
\sum_{j=1}^N\Re \vk_j\le dC_*^2\|V\|_p^2,\qqq if \qqq \Re  V\ge 0.
\end{aligned}
\]

\end{theorem}



\subsection{Real potentials}
 Note that some of the results stated in above theorems are new even
 for real-valued  potentials.
We consider Schr\"odinger operators $H=\D+V$ with real potentials
$V$ under the condition \er{V}. In this case all eigenvalues $\l_j$
and the numbers $z_j, \vk(\l_j)$ for all $ j=1,...,N$ are real. Thus
we have the same modified determinant $D(z)$ and Theorems
\ref{T1}-\ref{T3} hold true. Then from these results we obtain trace
formulae for real potentials.

\begin{corollary} \lb{T5} {\bf (Trace formulae.)} Let a potential $V$ be real and  satisfy
\er{V}. Then
\[
\lb{t1r} {\s(\T)\/2\pi}-B_0={1\/2\pi}\int_{-\pi}^{\pi}\log
|D(e^{it})|dt\ge 0,
\]
\[
\lb{zxjr} \sum_{j=1}^N\vk_j={d\/2\pi}\int_\T e^{-it}d\m(t),
\]
\[
\lb{zxr}
\begin{aligned}
  \sum_{j=1}^N\vk_{j}\l_{j}={1\/2}\Tr \,
V^2+{d^2\/4\pi}\int_\T e^{-2it}d\m(t),\\
\end{aligned}
\]
where $B_0=\log B(0)=\log \left(\prod_{j=1}^N |z_j|\right)<0$ and
$\vk_j=|\l_j^2-d^2|^{1\/2}\sign \l_j$.

\end{corollary}

From  Corollary \ref{T5} and Theorem \ref{T4} we obtain

\begin{corollary} \lb{T6} Let a potential $V$ be real and  satisfy
\er{V}. Then
\[
\lb{t1re} {\s(\T)\/2\pi}-B_0\le {C_*^2\/2}\|V\|_p^2,
\]
\[
\lb{zxjre} |\sum_{j=1}^N\vk_j|\le d C_*^2\|V\|_p^2,
\]
\[
\lb{zxre}
\begin{aligned}
  \sum_{j=1}^N\vk_{j}\l_{j}\le {1\/2}\Tr \,
V^2+{d^2\/4\pi}C_*^2\|V\|_p^2,\\
\end{aligned}
\]
where $B_0=\log B(0)=\log \left(\prod_{j=1}^N |z_j|\right)<0$ and
$\vk_j=|\l_j^2-d^2|^{1\/2}\sign \l_j$.

\end{corollary}

There are a lot of papers about eigenvalues of Schr\"odinger
operators in $\R^d$ with complex-valued potentials decaying at
infinity. Bounds  on sums of powers of eigenvalues were obtained in
\cite{FLLS06,LS09,DHK09,FLS16, F3, S10} and see references therein.
These bounds generalise the Lieb--Thirring bounds \cite{LT76} to the
non-selfadjoint case.

\medskip

For discrete Schr\"odinger operators most of the results were
obtained for the $\Z^1$ self-adjoint case, see, for example,
\cite{T89}. For  the nonself-adjoint case we mention \cite{BGK09}
and see references therein. Schr\"odinger operators with decreasing
potentials on the lattice $\Z^d, d>1$ have been considered by Boutet
de Monvel-Sahbani \cite{BS99}, Isozaki-Korotyaev \cite{IK12},
Isozaki-Morioka \cite{IM14}, Kopylova \cite{Ko10}, Korotyaev- Moller
\cite{KM17}, Rosenblum-Solomjak \cite{RS09}, Shaban-Vainberg
\cite{SV01} and see references therein.  Scattering on other graphs
was discussed by Ando \cite{A12}, Korotyaev-Saburova \cite{KS15} and
Korotyaev-Moller-Rasmussen \cite{KMR17}, Parra-Richard \cite{PR17}.

\medskip


\section {Preliminaries}
\setcounter{equation}{0}

\subsection{Trace class operators}
We recall some well-known facts.\\
$\bu$
Let $A, B\in \cB$ and $AB, BA\in \cB_1$. Then
\[
\label{2.1} {\rm Tr}\, AB={\rm Tr}\, BA,
\]
\[
\label{2.2} \det (I+ AB)=\det (I+BA).
\]
\\
$\bu$ Let  an operator-valued function $\O :\mD\to \cB_1$ be
analytic for some domain $\mD\ss\C$ and $(I+\O (z))^{-1}\in \cB$ for
any $z\in \mD$. Then for the function $F(z)=\det (I+\O (z))$ we have
\[
\label{S2F'z}
 F'(z)= F(z)\rm Tr \Omega (z)^{-1}\Omega '(z).
\]
$\bu$ In the case $A\in \cB_2$ we define the  modified  determinant
$\det_2 (I+ A)$ by
\[
\label{DA2} \det_2 (I+ A) =\det \rt((I+ A)e^{-A}\rt).
\]
The  modified  determinant satisfies (see (2.2) in Chapter IV,
\cite{GK69})
\[
\label{DA2x} |\det_2 (I+ A)|\le e^{{1\/2}\|A\|_{\cB_2}^2},
\]
and  $I+ A$ is invertible if and only if $\det_2 (I+ A)\ne 0$.

\subsection{Fredholm determinant}

Consider the bounded operators $V\in \cB_2$ and  $H_0$ acting in the
Hilbert space $\mH$. Define the operator $H=H_0+V$. Introduce the
resolvents
$$
R_0(\l)=(H_0-\l)^{-1},\qq \l\notin \s(H_0) \qqq {\rm and} \qqq
R(\l)=(H-\l)^{-1},\qq \l\notin \s(H).
$$
For  $V\in \cB_2$ we define the regularized determinant $\cD$ by
\[
\lb{dD2} \cD(\l)=\det \rt[(I+VR_0(\l))e^{-VR_0(\l)}\rt], \qq
\l\notin \s(H_0).
\]

Note that for any bounded operator $H$ and for large $\l$ we have
\[
\label{asR}
R(\l)=-{1\/\l}\sum_{n\ge 0}\rt({H\/\l}\rt)^n,\qqq |\l|>\|H\|,
\]
where the series is absolutely convergent.


\begin{lemma}
\label{TPD} Let operators $V\in \cB_2$ and $H_0\in \cB$ and the
modified determinant $\cD(\lambda)$ be defined by \er{dD2}. Then
$\cD(\lambda)$ is analytic in $\{\l\in \C:|\lambda|
> r_0\}$ for $r_0=\|H_0\|$. Moreover
\begin{equation}
\label{PD1} \cD(\lambda)=1+O(1/\l^2) \quad as \quad |\l|\to
{\infty},
\end{equation}
\begin{equation}
\label{PD3} \log \cD(\lambda) = -
\sum_{n=2}^{\infty}\frac{(-1)^n}{n}{\rm
Tr}\,\left(VR_0(\lambda)\right)^n,
\end{equation}
and
\begin{equation}
\label{PD4}
\begin{aligned}
\log \cD(\lambda) =-\sum _{n \geq 2}\frac{d_n}{\lambda^n}=
-{d_2\/\l^2}-{d_3\/\l^3}-{d_4\/\l^4}-....,\\
d_2={1\/2}\Tr\, V^2, \quad  d_n={1\/n}{\rm Tr}\,\big(H^n -
H_0^n-nH_0^{n-1}V\big), \quad n\ge 2,
 \end{aligned}
\end{equation}
where the right-hand side is uniformly convergent on $\{\l\in
\C:|\lambda|\ge r\}$ for $r=\|V\|+\|H_0\|+1$. In particular,
\begin{equation}
\label{PD5} d_3={1\/3}{\rm Tr}\,(3V^2H_0+V^3),\qq d_4={1\/4}{\rm
Tr}\,(2VH_0VH_0+4V^2H_0^2+4V^3H_0+V^4).
\end{equation}

\no ii) The function $\p(z)= \log \cD(\l(z))$ is analytic in
$\dD_{r}$ for some $r>0$ and has the following Tailor series
\[
\label{PD6} \p(z)=-\p_2z^2-\p_3z^3-\p_4z^4 +......, \qqq \as \qq
|z|<r,
\]
 and
\[
\label{PD7} \p_2=a^2d_2, \qqq \p_3=a^3d_3, \qqq \p_4=a^4
d_4-\p_2,....,
\]
here $a={2\/d}$ and the coefficients $d_n$ are given by \er{PD4}.
\end{lemma}
{\bf Remark.} Due to \er{PD1} we take the branch of $\log \cD$ so
that $\log \cD(\l )=o(1)$ as $|\l |\to {\infty}$.

{\bf Proof}. i)  We have
\begin{equation}
\label{PB1} \|\left(VR_0(\lambda)\right)^2\|_{\cB_1} \leq
{\|V\|_{\cB_2}^2\/|\lambda|^2} \quad {\rm for}\quad |\lambda|
> 2r.
\end{equation}
The Taylor series for the entire function  $e^{-T}$ and the estimate
\er{PB1} give at $T=VR_0(\l)$
$$
[(I+T)e^{-T}]=(I+T)(1-T+T^2O(1))=1-T^2+T^2O(1)=I+T^2O(1).
$$
Take $r_1 > 0$ large enough. Then  for $|\lambda| > r$, we have by
the resolvent equation
\begin{equation}
\label{RR0}
 R(\lambda) =R_0(\lambda)+
\sum_{n=1}^{\infty}(-1)^nR_0(\lambda) \rt(VR_0(\lambda)\rt)^{n}=
\sum_{n=0}^{\infty}(-1)^nR_0(\lambda) \rt(VR_0(\lambda)\rt)^{n},
\end{equation}
where the right-hand side is uniformly convergent on $\{\l\in
\C:|\lambda|\ge r\}$. By (\ref{S2F'z}), (\ref{PB1}) and using
\er{2.1}, we have for $|\l|> r$ the following
\begin{equation}
\lb{D16}
\begin{aligned}
&\cD'(\lambda)=-\cD(\lambda)\Tr \rt(Y(\lambda)Y_0'(\lambda))\rt)
=-\cD(\lambda)\Tr \rt(VR(\lambda)VR_0^2(\lambda)\rt)\\
&=-\cD(\lambda)\Tr \rt(R_0(\lambda)VR(\l)VR_0(\l)\rt)
=-\cD(\lambda)\Tr \rt(R(\l)(VR_0(\l))^2\rt).
\end{aligned}
\end{equation}
Thus (\ref{RR0}) gives
\begin{equation}
\begin{aligned}
\lb{deD} &(\log \cD(\lambda))'=-{\rm Tr}\,
\sum_{n=0}^{\infty}(-1)^nR_0(\lambda) \rt(VR_0(\lambda)\rt)^{n+2} =
-{\rm Tr}\, \sum_{n=2}^{\infty}(-1)^nR_0(\lambda)
\rt(VR_0(\lambda)\rt)^{n}.
\end{aligned}
\end{equation}
Then integrating and using
$$
 {d\/d\lambda}\rt(\Tr \rt(VR_0(\lambda)\rt)^{n}\rt)=n
 \Tr\, R_0(\lambda) \rt(VR_0(\lambda)\rt)^{n}
$$
we obtain \er{PD3}. The identities \er{D16} and $R=R_0-RVR_0$ imply
\begin{equation}
\begin{aligned}
(\log \cD(\lambda))' = -{\rm Tr}\,\Big(R(\lambda) -
R_0(\lambda)+R_0(\lambda)VR_0(\lambda)\Big) =-{\rm
Tr}\,\Big(R(\lambda) - R_0(\lambda)+VR_0^2(\lambda)\Big).
\end{aligned}
\nonumber
\end{equation}
Using the identity \er{asR}  we obtain
\begin{equation}
(\log \cD(\lambda))'= \sum_{n=0}^{\infty}{{\rm Tr}\,\big(H^n -
H_0^n-nH_0^{n-1}V\big)\/\lambda^{n+1}}=\sum_{n=2}^{\infty}{nd_n\/\lambda^{n+1}}.
\nonumber
\end{equation}
In view of \er{PD1}, we get \er{PD4}, \er{PD5}.

ii) Using \er{PD4}, \er{PD5} and the identity $\l={d\/2}(z+{1\/z})$
we obtain \er{PD6}, \er{PD7}. \BBox

{\bf Remark.} Consider the case when an operator $V\in \B_2$ and an
operator $H_0\in \B$ is self-adjoint. Then each zero of $\cD(\l)$
outside $\sigma(H_0)$ is eigenvalue of $H=H_0+V$ and its
multiplicity is a multiplicity of this eigenvalue.


For $\lambda \not\in \sigma(H_0)$, the eigenvalue problem
$(H-\lambda)u = 0$ is equivalent to $(I + (H_0-\lambda)^{-1}V)u =
0$, which has a non-trivial solution if and only if $\cD(\l)= 0$.

\

We will determine the asymptotics of $\log B(z)$ as $z\to 0$. For a
sufficiently small $z$ and for $t=z_j\in \dD$ for some $j$ we have
the following identity:
$$
\log {|t|\/t}{t-z\/1-\ol t z}=\log |t|+ \log \rt(1-{z\/t}\rt)-\log
(1-\ol t z)
=\log |t|-\sum_{n\ge1}\rt({1\/t^n}-\ol t^n\rt){z^n\/n}.
$$
Besides,
$$
\begin{aligned}
&  |1-|t|^n|\le n|1-|t||,\\
 &\big|t^{-n}-\ol t^n \big|\le \big|1-t^n| +\big|1-t^{-n} \big|\le
|1-t^n|\rt(1+{1\/|t|^n} \rt)\le
 |1-|t|^n|{2\/r_0^n}\le  |1-|t||{2n\/r_0^n},
\end{aligned}
$$
where $r_0=\inf |z_j|>0$.
 This yields
\[
\lb{asB}
\begin{aligned}
& \log  B(z)=\sum_{j=1}^N\log  {|z_j|\/z_j}{z_j-z\/1-\ol z_j z}
=\sum_{j=1}^N\rt( \log |z_j|+ \log \big(1-(z/z_j)\big)-\log (1-\ol
z_j
z)\rt)\\
&=\sum_{j=1}^N\log
|z_j|-\sum_{n=1}^\iy\sum_{j=1}^N\rt({1\/z_j^n}-\ol
z_j^n \rt){z^n\/n}=\log B(0)-b(z),\\
& b(z)=\sum_{n=1}^N\sum_{j=1}^\iy\rt({1\/z_j^n}-\ol z_j^n
\rt){z^n\/n}=\sum_{n=1}^N z^nB_n,\qqq
B_n={1\/n}\sum_{j=1}^N\rt({1\/z_j^n}-\ol z_j^n \rt),
\end{aligned}
\]
where  the function $b$ is analytic in the disk $\{|z|<{r_0\/2}\}$
and $B_n$ satisfy
$$
\begin{aligned}
 |B_n|\le {1\/n}\sum_{j=1}^N\rt|{1\/z_j^n}-\ol z_j^n \rt| \le
{2\/r_0^n}\sum_{j=1}^N  |1-|z_j||={2\/r_0^n}\cZ_D,
 \end{aligned}
$$
where $\cZ_D=\sum_{j=1}^\iy(1-|z_j|)$.
 Thus
$$
|b(z)|\le \sum_{n=1}^\iy |B_n|{|z|^n}\le
2\cZ_D\sum_{n=1}^\iy{|z|^n\/r_0^n}={2\cZ_D\/1-{|z|\/r_0}}.
$$


\section {Complex potentials}
\setcounter{equation}{0}

\subsection {Momentum representation}

Define the Fourier transformation $\F: \ell^2(\Z^d) \to L^2(\T^d)$
by
\[
f \to \hat f(k)=(\F f)(k)={1\/(2\pi)^{{d\/2}}}\sum_{n\in \Z^d}
f_ne^{i(n,k)},\qqq k=(k_j)_1^d\in \T^d,
\]
where $\T^d=\R^d/(2\pi \Z)^d$.  We need the so-called momentum
representation of the operator $H$:
\[
\wh H=\F H \F^*=\wh \D+\cV, \qqq \wh \D=\F \D \F^*, \qqq \cV=\F V
\F^*,
\]
\[
(\wh \D f)(k)=\wh \D(k)f(k),\qqq \wh \D (k)=\sum_1^d \cos k_j
\]
\[
(\cV f)(k)={1\/(2\pi)^{d\/2}}\int_{\T^d} \wh V(k-k')\wh f(k')dk',\qq
k=(k_j)_1^d\in \T^d,
\]
\[
\wh V(k)={1\/(2\pi)^{d\/2}}\sum_{n\in \Z^d}V_ne^{i(n,k)},\qqq
V_n={1\/(2\pi)^{d\/2}}\int_{\T^d} \wh V(k)e^{-i(n,k)}dk.
\]

\subsection {Preliminaries }

We define the Hardy space in the upper half-plane. Let $F(\l),
\l=\m+i\nu\in \C_+$ be analytic on $\C_+$. For $0<p\le \iy$ we say
$F$ belongs the Hardy space $ \mH_p=\mH_p(\C_+)$ if $F$ satisfies
$\|F\|_{\mH_p}<\iy$, where $\|F\|_{\mH_p}$ is given by
$$
\|F\|_{\mH_p}=\ca
\sup_{\nu>0}{1\/2\pi}\rt(\int_\R|F(\m+i\nu))|^pd\m\rt)^{1\/p} &
if  \qqq 0< p<\iy\\
 \sup_{\l\in \C_+}|F(\l)| & if \qqq p=\iy\ac .
$$
Note that the definition of the Hardy space $\mH_p$ involves all
$\nu=\Im \l>0$.

Recall that we have defined the new spectral variable $z\in \dD$ by
$$
\l=\l(z)={d\/2}\rt(z+{1\/z}\rt)\in \L=\C\sm [-d,d] ,\qqq z\in \dD.
$$
 The function $\l(z)$ has  the following properties:

{\it $\bu$  The function $\l(z)$ is a conformal mapping from $\dD$
onto the spectral domain $\L$.

$\bu$  $\l(\dD)=\L=\C\sm [-d,d] $  and   $\l(\dD\cap \C_\mp)=\C_\pm
$.

$\bu$  $\L$ is the cut domain with the cut $[-d,d]$, having the
upper side $[-d,d]+i0$ and the lower side $[-d,d]-i0$. The function
$\l(z)$ maps the boundary: the upper semi-circle onto the lower side
$[-d,d]-i0$ and the lower semi-circle onto the upper side
$[-d,d]+i0$.

$\bu$  The function $\l(z)$ maps the point $z=0$ to the point
$\l=\iy$.

$\bu$  The inverse mapping $z(\cdot ): \L\to \dD$ is given by
$$
\begin{aligned}
z={1\/d}\rt(\l-\sqrt{\l^2-d^2}\rt),\qqq \l\in \L,\\
z={d\/2\l}+{O(1)\/\l^3}\qqq as \qq |\l|\to \iy.
\end{aligned}
$$
}

\subsection {Proof of main theorems} We consider a Schr\"odinger operator $H=\D+V$
on $\ell^2(\Z^d)$. We
assume that the potential $V$ is  complex  and  satisfies the
condition \er{V}. We present preliminary results. Recall that we can
rewrite the modified determinant $\cD(\l), \l\in \L$ in the form
$$
\cD(\l)=\det \rt[(I+Y_0(\l)e^{-Y_0(\l)}\rt], \qqq z\in \L,
$$
where the operator $Y_0(\l)$ is given by
\[
\label{DeY}
\begin{aligned}
Y_0(\l)=|V|^{1\/2} R_0(\l) V^{1\/2}, \qq \
V^{1\/2}=|V|^{1\/2}e^{i\arg V}, \qq \l\in\L=\C \sm \s(H_0).
\end{aligned}
\]

{\bf Proof of Theorem \ref{T1}.} Recall that the modified
determinant $D(z)=\cD(\l(z)), z\in \dD$. The determinant $D(z), z\in
\dD$ is well defined, since $V\in \cB_2$. It is well known that if
$\l_0\in \L$ is an eigenvalue of $H$, then $z_0=z(\l_0)\in \dD$ is a
zero of $D$ with the same multiplicity. We recall needed  results
from \cite{KM17}:

{\it Let the potential $V$ satisfy \er{V}. Then the operator-valued
function $Y_0: \C\sm [-d,d]\to \cB_2$ is analytic and H\"older
continuous up to the boundary. Moreover, it satisfies
\begin{equation}
\label{Y01}
\begin{aligned}
 \|Y_0(\l)\|_{\cB_2}\le C_*\|V\|_p,\qqq \forall \ \ \l\in
\L,
\end{aligned}
\end{equation}
\begin{equation}
\begin{aligned}
\label{Y02} \|Y_0(\l)-Y_0(\m)\|_{\cB_2}\le
C_\a|\l-\m|^\a\|V\|_p,\qqq \forall \ \ \qq \l,\m\in \ol \C_\pm,
\end{aligned}
\end{equation}
 where $C_\a$ is some constant and the constant $C_*$ is  defined by
\[
\label{Cpd}
\begin{aligned}
C_*=C_{p,d}+C_d^0 \G(p,d),\qq C_{1,d}=1,\qq C_{p,d}=p^{d(p-1)\/2 p},\\
\G(p,d)=\big(3+2\vk  \big)^{d(p-1)\/p},\qqq
C_d^0=\ca 16\\
           4\\
           {14\cdot 2^{d\/4}\/d-4}\ac,\qqq
 \vk=\ca {6(p-1)\/6-5p}\ & if \ d=3\\
      \rt({5p-1\/4-3p}\rt)^{5p-4\/4(p-1)}\ & if \ d=4\\
      {3d(p-1)\/3d-(2d+1)p}\ & if \ d\ge 5\ac
\end{aligned}
\]

}

Due to results \er{Y01}-\er{Y02} the operator-valued function
$Y_0(\l): \C_\pm\to \cB_2$ is analytic in the upper half-plane
$\C_+$ and is H\"older up to the boundary. Then the determinant
$\cD(\l)$ is analytic in the upper half-plane $\C_\pm$ and H\"older
up to the boundary, and satisfies
\[
\label{cD1}
\begin{aligned}
\|\cD\|_{\mH_\iy(\C_\pm)}\le e^{C_*^2\|V\|_{p}^2/2},
\end{aligned}
\]
where the constant $C_*$ is defined in \er{Cpd}. The function
$\cD(\l)$ has asymptotics \er{PD1}, then all zeros of $\cD(\l)$  are
uniformly bounded, which yields $\sum |\Im \l_j|<\iy$.

Consider the function $f(\vk)=\cD(\l(\vk)), \vk\in \K$, where
$\l(\vk)=\sqrt{\vk^2+d^2}$ is the conformal mapping $\K\to \L$. The
function $f\in \mH_\iy(\K_\pm)$, where $\K_\pm=\{\pm \Re z>0\}$.
Repeating arguments for the function $\cD(\l)$  we obtain $\sum |\Re
\vk(\l_j)|<\iy$.

Thus similar arguments give that the operator-valued function
$Y_0(\l(z)): \dD\to \cB_2$ is analytic in the unit disc $\dD$ and is
H\"older up to the boundary. Then   the determinant $D(z)$ is
analytic in the unit disc $\dD$ and H\"older up to the boundary, and
satisfies \er{D1}.

Furthermore, due to  Lemma \ref{TPD}    the function $\p(z)=\log
D(z)$ defined by $\log D(0)=0$ is analytic in the disc $\dD_{r_0}$
with the radius $r_0>0$ defined by $r_0=\inf |z_j|>0$ and has the
Tailor serious as $|z|<r_0$ given by \er{D3}.
 \BBox

We now consider  the canonical representation \er{cfD} (see,
\cite{Koo98}, p. 76):

  {\it Let a function $f\in \mH_p, p\ge 1$ and let $B$ be its
Blaschke product. Then $f$ has the form
\[
\lb{cf}
\begin{aligned}
 f(z)=B(z)e^{ic-K_\s (z)}e^{K_f(z)},\\
K_\s(z)={1\/2\pi}\int_{-\pi}^{\pi}{e^{it}+z\/e^{it}-z}d\s(t),\\
K_f(z)= {1\/2\pi}\int_{-\pi}^{\pi}{e^{it}+z\/e^{it}-z}\log
|f(e^{it})|dt,
 \end{aligned}
\]
for all $|z|<1$, where $c$ is real constant and $\log |f(e^{it})|\in
L^1(-\pi,\pi)$ and $\s=\s_f\ge 0$ is a singular measure on
$[-\pi,\pi]$ such that $\supp \s\ss \{t\in [-\pi,\pi]:
D(e^{it})=0\}$.}


\noindent We define the inner function and the outer function (after
Beurling) in the disc by
$$
\begin{aligned}
& f_{in}(z)=B(z)e^{ic-K_\s (z)} \qq & the \ inner \ factor\ of\ f,\\
&f_{out}(z)=e^{K_f(z)} \qq & the \ outer \ factor\ of\ f,\\
 \end{aligned}
$$
for $|z|<1$. Note that we have  $| f_{in}(z)|\le 1$, since $d\s\ge
0.$

We describe the canonical representation of the determinant
$D(z)=\cD(\l(z)), z\in \dD$.

{\bf Proof of Corollary \ref{T2}.} Theorem \ref{T1} implies $D\in
\mH_\iy$. Therefore the  canonical representation \er{cf} gives
\[
\lb{Dx1} D(z)=B(z)e^{ic-K_\s (z)}e^{K_D(z)},\qqq z\in \dD.
\]
In order to prove \er{cfD} we need to show $e^{ic}=1$. From \er{Dx1}
at $z=0$ we obtain
$$
1=D(0)=B(0)e^{ic-K_\s (0)}e^{K_D(0)}.
$$
Since $B(0), K_\s (0), K_f(0)$ and $c$ are real  we obtain
$e^{ic}=1$. \BBox

We describe trace formulae.

{\bf Proof of Theorem \ref{T3}. (The trace formulae.)} Due to the
canonical representation \er{cfD},  the function
$D_B(z)={D(z)\/B(z)} $ has no zeros in the disc $\dD$ and satisfies
\[
\lb{Si} \log
D_B(z)={1\/2\pi}\int_{-\pi}^{\pi}{e^{it}+z\/e^{it}-z}d\m(t),\qq
z\in\dD,
\]
where the measure $d\m=\log |f(e^{it})|dt-d\s(t)$. In order to show
\er{t1}--\er{t3} we need the asymptotics of the Schwatz integral
$\log D_B(z)$  as $z\to 0$. The following identity holds true
\[
\lb{ts1}
{e^{it}+z\/e^{it}-z}=1+{2ze^{-it}\/1-ze^{-it}}=1+2\sum_{n\ge 1}
\big({ze^{-it}}\big)^n=
1+2\big({ze^{-it}}\big)+2\big({ze^{-it}}\big)^2+.....
\]
for all $(t,z)\in \pa \dD\ts \dD$. Thus \er{Si}, \er{ts1} yield the
Taylor series at $z=0$:
\[
\lb{asm}
{1\/2\pi}\int_{-\pi}^{\pi}{e^{it}+z\/e^{it}-z}d\m(t)={\m(\T)\/2\pi}+\m_1z+
\m_2z^2+\m_3z^3+\m_4z^4+...\qqq as \qqq |z|<1,
\]
where
$$
\m(\T)=\int_0^{2\pi}d\m(t),\qqq \m_n={1\/\pi}\int_0^{2\pi}e^{-in\vt}
d\m(t),\qqq n\in \Z.
$$

We have the identity $ \log D(z)=\log B(z)+
{1\/2\pi}\int_{-\pi}^{\pi}{e^{it}+z\/e^{it}-z}d\m(t)$ for all $z\in
\dD_{r_0}$. Combining asymptotics \er{D3}, \er{asB} and \er{asm} we
obtain \er{t1}-\er{t3}. In particular, we have \er{t4} and $-\log
B(0)={\m(\T)\/2\pi}={1\/2\pi}\int_0^{2\pi}\log
|f(e^{it})|dt-{\s(\T)\/2\pi}\ge 0$.

Recall that $\vk=\sqrt{\l^2-d^2}$. We have the following identities
for $z\in \dD$ and $\l\in \L$:
\[
\label{Et1}
\begin{aligned}
{2}\l=d(z+{1\/z}),\qqq dz={\l-\vk},\qqq d(z-{1\/z})=-{2}\vk.
\end{aligned}
\]
These identities yield
\[
\label{Et2}
\begin{aligned}
 d\rt({1\/z}-\ol z\rt)={2}\l-2d\Re z,\qqq {d\/2}\Im\rt({1\/z}-\ol z\rt)=\Im \l,\\
d\rt({1\/z}-\ol z\rt)=2\vk +2id\Im z,\qqq {d\/2}\Re \rt({1\/z}-\ol z\rt)=\Re \vk,\\
 {d\/2}\rt({1\/z}-\ol z\rt)= \Re \vk+i \Im \l.
 \end{aligned}
\]
Let $\vk_j=\sqrt{\l_j^2-d^2}$. Then from \er{t2} we get
$$
 B_1=\sum_{j=1}^N\rt({1\/z_j}-\ol z_j
\rt)={1\/\pi}\int_\T e^{-it}d\m(t),
$$

$$
\begin{aligned}
{d\/2\pi}\int_\T e^{-it}d\m(t)=\sum_{j=1}^N{d\/2}\rt({1\/z_j}-\ol
z_j\rt)
=\sum_{j=1}^N\rt(\Re \vk_j+i \Im \l_j\rt)
 \end{aligned}
$$
and thus
$$
\sum_{j=1}^N \Re\sqrt{\l_j^2-d^2}={d\/2\pi}\int_\T \cos
t\,d\m(t),\qqq \sum_{j=1}^N \Im\l_j=-{d\/2\pi}\int_\T  \sin t\,d
\m(t),
$$
We show \er{zx}.  The function $\l=\l(z)={d\/2}\rt(z+{1\/z}\rt)$
satisfies
\[
\label{xx1}
\begin{aligned}
 d\rt({1\/z}+\ol z\rt)=d\rt({1\/z}+z+(\ol z-z)\rt)=2\l-2d\Im z,\qqq {d\/2}\Re \rt({1\/z}+\ol z\rt)=\Re \l,\\
d\rt({1\/z}+\ol zrt)=d\rt(z+2\vk+\ol z\rt)  =2\vk+2d\Re z,\qqq
{d\/2}\Im \rt({1\/z}+\ol
z\rt)=\Im \vk,\\
{d\/2}\rt({1\/z}+\ol z\rt)= \Re \l+i \Im \vk.
 \end{aligned}
\]
Let $\vk=\vk_1+i\vk_2$ and $\l=\l_1+i\l_2$.  Then we obtain
\[
\begin{aligned}
{d^2\/4}\rt({1\/z^2}-\ol z^2\rt)={d^2\/4}\rt({1\/z}-\ol
z\rt)\rt({1\/z}+\ol z\rt)
\\
=(\vk_1+i\l_2)(\l_1+i\vk_2)
=(\vk_1\l_1-\l_2\vk_2)+i(\l_1\l_2+\vk_1\vk_2).
\end{aligned}
\]
Note that $ \sign \vk_1\l_1=\sign \l_2\vk_2 $ and $\sign
\l_1\l_2=\sign \vk_1\vk_2. $ Thus we get
\[
\lb{t4x}
\begin{aligned}
{1\/2}\Tr \, V^2+{d^2\/4\pi}\int_\T
e^{-i2t}d\m(t)={d^2\/4}\sum_{j=1}^N\rt({1\/z_j^2}-\ol z_j^2 \rt)
\\
=
\sum_{j=1}^N\rt( (\vk_{j1}\l_{j1}-\l_{j2}\vk_{j2})+i(\l_{j1}\l_{j2}+\vk_{j1}\vk_{j2})  \rt),\\
\end{aligned}
\]
and then
\[
\lb{123}
\begin{aligned}
  \sum_{j=1}^N(\vk_{j1}\l_{j1}-\l_{j2}\vk_{j2})={1\/2}\Tr \,
\Re V^2+{d^2\/4\pi}\int_\T \cos 2td\m(t),\\
\sum_{j=1}^N(\l_{j1}\l_{j2}+\vk_{j1}\vk_{j2})={1\/2}\Tr \, \Im
V^2+{d^2\/4\pi}\int_\T \sin 2td\m(t).
\end{aligned}
\]

\BBox

\

 {\bf Proof of Theorem \ref{T4}. Estimates.}
 The simple inequality
$1-x\le -\log x$ for $\forall \ x\in (0,1]$, implies
$-B_0=-B(0)=-\sum \log |z_j|\ge \sum (1- |z_j|)$. Then substituting
the last estimate and the estimate \er{D1} into the first trace
formula \er{t1}  we obtain \er{eiV}.

We show \er{em}. Let for shortness $C=C_*^2\|V\|_p^2$.  From \er{D1}
and \er{eiV}  we obtain
\[
\lb{mnC}
\begin{aligned}
  {1\/2\pi}\big|\int_\T
e^{-int}d\m(t)\big|\le {1\/2\pi}\int_\T
\rt({C\/2}dt+d\s\rt)={C\/2}+{\s(\T)\/2\pi}\le C, \qqq \forall \ n\in
\Z.
\end{aligned}
\]
In order to determine the next two estimates we use the trace
formula \er{zxj}. From \er{zxj} and \er{mnC} we obtain
\[
\lb{2es}
\begin{aligned}
 |\sum_{j=1}^N\big(\Re \vk_j+i \Im \l_j\big)|\le
{d\/2\pi} \big|\int_\T e^{-int}d\m(t)\big|
 =d  C_*^2\|V\|_p^2.
\end{aligned}
\]
If  $\Im V\ge 0$ (or $\Re V\ge 0$), then $\Im \l_j\ge 0$ ($\Re
\l_j\ge 0$)  and the estimates \er{2es} gives \er{eiV2}.
 \BBox


%

\footnotesize \no \noindent \textbf{Acknowledgments.}  EK  is also
grateful to A. Alexandrov (St. Petersburg) and K. Dyakonov
(Barcelona) for  useful comments about Hardy spaces. Our study was
supported by the RSF grant No 15-11-30007.


\begin{thebibliography}{99}


\bibitem[A12]{A12} K. Ando, Inverse scattering theory for
discrete Schr\"odinger operators on the hexagonal lattice, Ann.
Henri Poincar\'e, 14 (2013), 347--383.


\bibitem[BGK09]{BGK09} A. Borichev, L. Golinskii, S. Kupin,
A Blaschke-type condition and its application to complex Jacobi
matrices. Bull. London Math. Soc., 41 (2009), 117--123.







\bibitem[BS99]{BS99}
A. Boutet de Monvel; J. Sahbani, On the spectral properties of
discrete Schr{\"o}dinger operators : (The multi-dimensional case),
Review in Math. Phys., 11 (1999), 1061-1078.


\bibitem[DHK09]{DHK09} M. Demuth; M. Hansmann; G. Katriel,
On the discrete spectrum of non-selfadjoint operators.  J. Funct.
Anal. 257 (2009), no. 9, 2742--2759.


\bibitem[F3]{F3} R. L. Frank, Eigenvalue bounds for Schr\"odinger
operators with complex potentials. III. Preprint (2015),
http://arxiv.org/pdf/1510.03411v1.pdf

\bibitem[FLS16]{FLS16} R. L. Frank, A. Laptev, O. Safronov,
\textit{On the number of eigenvalues of Schr\"odinger operators with
complex potentials}, J. London Math. Soc. (2016) doi:
10.1112/jlms/jdw039.

\bibitem[FLLS06]{FLLS06} R. L. Frank, A. Laptev, E. H. Lieb, R. Seiringer,
Lieb--Thirring inequalities for Schr\"odinger operators with
complex-valued potentials. Lett. Math. Phys. 77 (2006), 309--316.


\bibitem[FS14]{FS14} R. L. Frank, J. Sabin, Restriction theorems for orthonormal functions,
Strichartz inequalities and uniform Sobolev estimates. Preprint
(2014), http://arxiv.org/pdf/1404.2817.pdf



\bibitem[GK69]{GK69} I. Gohberg; M. Krein,
Introduction to the theory of linear nonselfadjoint operators.
Translated from the Russian, Translations of Mathematical
Monographs, Vol. 18 AMS, Providence, R.I. 1969.

\bibitem[IK12]{IK12} H. Isozaki; E. Korotyaev, Inverse Problems,
Trace formulae for discrete Schr\"odinger Operators,  Annales Henri
Poincare, {\bf 13} (2012), No 4 ,  751--788.



\bibitem[IM14]{IM14} H. Isozaki; H. Morioka, A Rellich type theorem for
discrete Schr\"odinger operators,   Inverse Probl. Imaging, {\bf 8}
(2014), no. 2, 475–-489.


\bibitem[Koo98]{Koo98} P. Koosis, Introduction to $H_p$ spaces, 115
Cambridge Tracts in Mathematic, 1998.

\bibitem[Ko10]{Ko10} E.A. Kopylova,  Dispersive estimates
 for discrete Schr\"odinger and Klein-Gordon equations,
 St. Petersburg Math. J.,  21
(2010), no. 5, 743–-760.


\bibitem[KL16]{KL16} E. Korotyaev; Laptev, A.
Trace formulae for Schr\"odinger operators with complex-valued
potentials on cubic lattices, preprint 2016.

\bibitem[KM17]{KM17} E. Korotyaev; J.S. Moller,
Weighted estimates for the Laplacian on the cubic lattice, preprint
2017.

\bibitem[KMR17]{KMR17} E. Korotyaev; J.S. Moller; M.G. Rasmussen,
Estimates for metric Laplacians on square  lattice,  2017

\bibitem[KS15]{KS15}   E. Korotyaev, N. Saburova,
Scattering on periodic metric graphs,   arXiv:1507.06441.


\bibitem[LS09]{LS09} A. Laptev, O. Safronov, Eigenvalue estimates
for Schr\"odinger operators with complex potentials. Comm. Math.
Phys. 292 (2009), 29--54.




\bibitem[LL76]{LL76}
E. Lieb; M. Loss, Analysis, AMS, Graduete Studies in Math., 14
(1976).

\bibitem[LT76]{LT76} E. H. Lieb, W. Thirring,
Inequalities for the moments of the eigenvalues of the Schr\"odinger
Hamiltonian and their relation to Sobolev inequalities. Studies in
Mathematical Physics. Princeton University Press, Princeton (1976),
pp. 269--303.

\bibitem[PR17]{PR17}
D. Parra; S. Richard, Spectral and scattering theory for Schrodinger
operators on perturbed topological crystals. arXiv:1607.03573.


\bibitem[RS09]{RS09} G. Rosenblum, M. Solomjak, On the spectral
estimates for the Schr{\"o}dinger operator on ${\Z}^d$, $d \geq 3$,
Problems in Mathematical Analysis, No. 41, J. Math. Sci. N. Y. 159
(2009), No. 2, 241--263.


\bibitem [S10]{S10}  O.  Safronov,  On a sum rule for Schr\"odinger
 operators with complex potentials,  Proc. Amer. Math. Soc. 138 (2010),
 no. 6, 2107--2112.


\bibitem[SK05]{SK05}     A. Stefanov; P.G. Kevrekidis, Asymptotic behaviour
of small solutions for the discrete nonlinear Schr\"oodinger and
Klein–Gordon equations, Nonlinearity, 18 (2005), pp. 1841–-1857.

 \bibitem[SV01]{SV01}  W. Shaban, B. Vainberg, Radiation conditions for
the difference Schr\"odinger operators, J. Appl. Anal., 80 (2001)
525--556.

\bibitem[T89]{T89} M. Toda, Theory of Nonlinear Lattices, 2nd. ed., Springer,
Berlin, 1989.



\end{thebibliography}
\end{document}